UDC 517.977.5+004.032.26 – Optimal control and artificial neural networks

# EXTENT OF ERROR CONTROL IN NEURAL NETWORKS

© 2016 A.V. Ignatenkov[1], A.M. Olshansky[2]

*Alexander Vladimirovich Ignatenkov, CEO for Information Technologies, LLC "Scientific Technological Center for Railway Management", postgraduate student of Rail Automation, Telemechanics and Communication Department of Samara State Transport University.*
*E-mail: a.ignatenkov@gmail.com*
*Alexey Mikhailovich Olshansky, PhD in Engineering,*
*doctoral student of Rail Automation, Telemechanics and Communication Department of Samara State Transport University E-mail: lexolshans@gmail.com*

The article sets and solves the task to control an error of the artificial neural network with variable signal conductivity.

Recently artificial neural networks (ANN) have been used to solve scientific and practical problems such as control synthesis for different systems, process identification, data handling, information recovery, process scheduling and creating images which is a very unordinary application. It holds for both networks' structure or topology and learning algorithms and methods. Nevertheless ANNs with their own behavior are kept in the same status of entirely experiential systems with no theoretical methods to apply. For example S. Haykin [8] notes that rigor research of ANN's behavior is a complex scientific problem.

Every ANN is a complex system with variety of characteristic features which allow us to consider ANNs as dynamic systems. There are only a few investigations where ANN's are considered as dynamic or optimal control systems.

As pioneering studies we should mention papers on the optimal neural network structure. In [5] the authors use an evolutionary algorithm of simultaneous adjustment of weight structure and links between neurons to optimize the ANN structure. Complexity of objective estimation of ANN topology is one of the suppositions of this kind of algorithms [5, pp.35-36]. The main idea of this algorithm is to encode the values of weights and the number of neurons in ANN genotype and to run crossover, mutation and selection operators. Populations of ANNs are being hybridized and selected with population reduction and elite selection. Meanwhile the author [5, p.103] notes that there is no optimal set of parameters of any algorithm to solve ANNs creating problem with the maximal rate of efficiency.

We agree with the observation made in the paper mentioned above that traditional learning methods and approaches fail to determine components of ANN output when we feed train samples and operate with output data sequence. It is suitable case for ANN with changeable signal propagation during scheduling problem operation. ANN with changeable signal propagation is able to analyze information about value of error from the last layer only. That's why it is urgent to find other approaches to control the ANN. Besides, no matter what ANN configuration we have, its task is to provide efficient learning.

An attempt to create schedules using Bellman equation [10] is analogous to the task mentioned above. It is implemented with stepwise schedule using any optimal part of a schedule when the schedule is generally optimal. Formulating the task we assume that:

$T_i$ – is the required service time,

$W_i$ is the time required to meet the *i*-request,

$a_i$ is the penalty for a waiting time unit after receiving the *i*-request, then waiting time for the result is $V_i = T_i + W_i$.

Waiting time penalty could be written as follows:

$$C = \sum_i a_i * (W_i + T_i)$$

Let's denote schedule completion time of every request as $D_i$. If we want the work to be done on time, we calculate penalties as maximal or mean upset of request *i* completion time. Performance criterion is defined as

$C = \max_i a_i \{\max[0, (V_i - D_i)]\}$,

where $\max[0, (V_i - D_i)]$ is the schedule time excess.

The schedule time may be exceeded both for particular tasks and most requests. In this case let's use General Penalty Criterion:

$$C = \sum_{i=1}^{n} a_i \{\max[0, (V_i - D_i)]\}$$

Then Bellman equation to calculate minimal amount of penalties for optimal service of first *i* requests and mean $T_i$ can be written as:

$$C_n(T) = \min_i [C_{n-1}(T - T_i) + a_i T] \qquad (1)$$

Expression (1) means that with the optimal sequence for all *n* tasks for the total *T* time the sequence of solving *n-1* tasks for the *T-$T_i$* time is to be optimal, too. Then the total minimum of penalties is a set when a permutation of every two tasks does not reduce the total number of penalties.

From (1) we may get

$$\frac{T_{i-1}}{a_{i-1}} \le \frac{T_i}{a_i} \qquad (2)$$

This way we got the rule (2) for composing an optimal schedule in order of increasing of penalty value/penalty for a portion of waste time rate. If waste penalties are equal we should compose schedules in order of $T_i$ values increasing. So we may prove that this approach guarantees us minimal average number of requirements in a system and minimal average waste time to requirements handling [11].

An exponential growth of number of possible permutations is a disadvantage of this approach. Besides this problem is an ideal hardly applicable to real railway traffic. In practice we have to consider both formalized and unformulated scheduling constraints. Special instruments are to be developed for building the schedule as a whole. ANN is a sample of such methods.

Another set of research contains papers which attempt to synthesize of automation control theory and ANN approach. In [3] optimal weights consequences per a time moments are constructed. A two point-boundary non-linear problem is solved resulting in learning rules for ANN.

Optimal ANN weight matrix is the sum of outer products between desired (n*1) and output vectors of the concrete subsystem. Every epoch has its own optimal weight matrix for ANN. In ideal case this matrix in final time moment looks like J.J. Hopfield symmetric matrix for associative memory [1].

Neuron evolution equation is

$$\frac{dX}{dt} = -X(t) + W(t)g(x(t)), \qquad (3)$$

where W(t) is a block-diagonal matrix, every block in the matrix is a W1(t) matrix,

g(x) is neuron output.

Initial conditions for (3) is concatenated input vector.

Performance functional is a correlation between neuron output and desired output in a final moment. It penalizes uncorrelatedness between neural activation function output and desired output.

Optimal control in this case can be solved as Lagrange problem to create program control. These results refer to the particular case of ANN described in [3]. Besides there are doubts about correctness of ANN phase vector velocity estimation as a difference between desired output and ANN output per epoch (3).

This approach doesn't operate with non-fully controllable ANNs while we have a group of ANNs which are output-controllable but not state-controllable.

Paper [4] deals with developing of intelligent control using so called dynamic neural networks. A dynamic neural network may be described by standard dynamic equations but in a standard scheme of a system we can observe dynamic chains such as an integration chain, an amplification chain, negative feedbacks. If characteristic matrices of this system are in concrete interval of values then ANN may demonstrate limited cycles, chaotic evolution and essential nonlinear evolution. In [4] ANNs are considered as self-organizing systems with quadratic performance functional, the authors of [4] also modify BFGS optimization algorithm. According to the [4] the growth of number of neurons is described by equation similar to (3) and almost heat-transfer equation signal distribution. If a neuron's axon get the environment of a next neuron and direction of axon's growth is congruent with mediator's concentration gradient the neuron starts expanding. But this approach doesn't describe functioning of ANN.

Paper [7] is a similar to author's field of research. A model in [7] is congruent with the model in [3] with constrained weights and ANN is a control system with delays. Both [3] and [7] are not completely suitable for our ANN topology with variable signal conduction (ANNV).

The papers on control synthesis of neural networks show that this subject needs further investigation. A very interesting and potentially novel problem is to evaluate ANNV [2] state based on observable output or on error signal level. Potential benefits of this problem are in composing schedules of various processes.

In this case we face with formal description of processes in ANNV in equation form or in matrix form. In the article authors present an attempt to synthetize optimal control strategy for ANNV as an open-loop system and then as a control system with feedback with quadratic performance criterion based on Bolza problem. In [9] the authors showed complete nonsatisfiability of maximal speed-in-action task and linear forms of performance criteria.

Let's consider a multilayer ANNV [2] described as (4):

$$\frac{dE}{dt} = \sum_{i=1}^{m}(a_i cos(w_i t) + b_i sin(w_i t)) + u(E,t) \quad (4)$$

where
m is the number of harmonics of error signal, m=7
$a_i$, $b_i$ – cosinus-coeffs and sinus-coeffs respectively,
$w_i$ – oscillation frequencies defined by spectral analysis of error signal,
$u(E,t)$ – some required control,
ANNV error signal could be represented as a set of harmonics with their own frequencies.

Assume that we may apply such control strategy $u$ to (4) to transform our ANNV from state 1 in $t_0$ to the desired state in time moment $t_1$. The control start moment $t_0$ is a moment when ANNV's trajectory E(t) is crossing a surface with preset level of error Δ. The final control moment $t_1$ is unfixed and it is value to define during the control problem's solving. A term of control strategy interrupting is E(t)< Δ.

*Program control synthetizing.*
Let's describe our network as follows:

$$\frac{dE}{dt} = \sum_{i=1}^{m}(a_i cos(w_i t) + b_i sin(w_i t)) + u(E,t) \quad (5)$$

All components in formula (5) are equivalent to (4).
Let us define the performance functional to minimize as

$$I = \int_{t_0} u^2(t)dt + E(t_1) \to min \quad (6)$$

For the problem (5)-(6) the Hamiltonian expression is written as:

$$H(E,t,u,\Psi) = \Psi_1(t) * \left(\left(\sum_{i=1}^{m}(a_i cos(w_i t) + b_i sin(w_i t))\right) + u(t)\right) - u^2(t) - E(t_1) \quad (7)$$

where $\Psi_1(t)$ is a special slave function. Assume that there are no constraints on control for problem (5)-(6).

According to the Pontryagin maxima principle the structure of optimal control strategy is

$$\frac{\partial}{\partial u}H(E,t,u,\Psi) = \Psi_1(t) - 2u(t) \quad (8)$$

Let us show that formula (8) is the maximum of Hamiltonian by control variable let's note (9):

$$\frac{\partial^2 H(E,t,u,\Psi)}{\partial u^2} = -2 < 0 \quad (9)$$

Find the structure of optimal control strategy solving the equation: (8)=0:

$$u_{opt}(t) = \frac{\Psi_1(t)}{2} \quad (10)$$

Let's formulate the conventional system of Pontryagin maxima principle based on the Hamiltonian in (7):

$$\begin{cases} \frac{\partial E}{\partial t} = \frac{\partial H}{\partial \Psi} = u(t) + \left(\sum_{i=1}^{m}(a_i\cos(w_i t) + b_i\sin(w_i t))\right) = \frac{\Psi_1(t)}{2} + \left(\sum_{i=1}^{m}(a_i\cos(w_i t) + b_i\sin(w_i t))\right), \\ E(t_0) = E_0, \\ \frac{\partial \Psi}{\partial t} = -\frac{\partial H}{\partial E} = -1 \end{cases}$$

(11)

The second equation in (11) contains initial conditions for equation in (5).

To get missing condition to solve (11) let's check transversality conditions. In general case transversality conditions are formulated as:

$$\delta E(t_1) - H(t_1) * \delta t_1 + \Psi_1(t_1) * \delta E = 0 \tag{12}$$

or

$$\delta E(t_1) - \left[\Psi_1(t) * \left(\left(\sum_{i=1}^{m}(a_i\cos(w_i t) + b_i\sin(w_i t))\right) + u(t)\right) - u^2(t) - E(t_1)\right] * \delta t_1 +$$
$$+\Psi_1(t_1) * \delta E = 0 \tag{13}$$

As there are no constraints for error value in the final control moment variations of error level $\delta E$ variations are arbitrary.

Choosing zero variations we have:

$$-\left[\Psi_1(t) * \left(\left(\sum_{i=1}^{m}(a_i\cos(w_i t) + b_i\sin(w_i t))\right) + u(t)\right) - u^2(t) - E(t_1)\right] * \delta t_1 = 0 \tag{14}$$

The right point of the trajectory is unfixed. It is equivalent to $\delta t_1 \neq 0$. Divide (14) by $\delta t_1$.

Let's substitute in formula (14) control $u(t)$ as it is in formula (10). We have

$$-\Psi_1(t) * \sum_{i=1}^{m}(a_i\cos(w_i t) + b_i\sin(w_i t)) + \frac{-\Psi_1(t)^2}{2} + \frac{\Psi_1(t)^3}{4} = E(t_1) \tag{15}$$

Solving the third equation in (11) we get

$$\Psi_1(t) = -t + C \tag{16}$$

In formula (16) for every final control moment we can calculate $C$ value using (15).

Then $u_{opt}(t) = \frac{-t+C}{2}$. Substituting u(t) in the first equation in (11) we have:

$$\frac{\partial E}{\partial t} = \frac{-t+C}{2} + \left(\sum_{i=1}^{m}(a_i\cos(w_i t) + b_i\sin(w_i t))\right) \tag{17}$$

Having integrated (17) we get

$$E(t) = \frac{Ct}{2} - \frac{t^2}{4} + \cdots + \frac{a_i}{w_i}\sin(w_i t) - \frac{b_i}{w_i}\cos(w_i t) + C_1 \tag{18}$$

where $C_1$ is calculated from the second equation in (11) for every known $C$ value.

The solution given in (18) has its own feature. If $t$ value is increasing we have the value $t$ which can be characterized with E(t) =0, but if $t$ value increases more and more E(t) is tending to minus infinity. This fact is in physical contradiction with sense of E(t). That's why we stop control strategy when E<E(t₁). Proceeding for this premise, let us solve the problem as an ANNV feedback control problem.

*Feedback control synthetizing.*

Let's consider a quadratic functional to minimize as (19):

$$I = \int_{t_0} u^2(E, W, t)dt + E(W, t_1, u) \to min \tag{19}$$

Every solution of ANN control problem contains feedback control strategy $u^*(E, W, t)$, final control time moment $t_1^*$ and an optimal error decreasing trajectory $E^*(t, W, u)$ under the optimal control strategy.

The optimal control strategy should be implemented algorithmically by the ANNV's structures and ANNV's behavior.

The form of boundary condition (20) is determined via suggestion about different error level for different final control time moments:

$$\Phi(t_1, E) = E \tag{20}$$

where $\Phi(t_1, E)$ is a boundary condition function.

For the problem (4) and criterion (5) with boundary condition (20) the Bellman equation looks like:

$$\Phi^Б = \frac{\partial \Phi}{\partial t} + \frac{\partial \Phi}{\partial E} * \left(\left(\sum_{i=1}^{m}(a_i cos(w_i t) + b_i sin(w_i t))\right) + u(E,t)\right) - u^2(E,t) \quad (21)$$

Let's find Bellman equation derivative and let's set it to 0 to find the optimal control structure:

$$\frac{\partial \Phi^Б}{\partial u} = \frac{\partial \Phi}{\partial E} - 2u(E,t) = 0 \quad (22)$$

The structure of optimal control is

$$u(E,t) = \frac{1}{2}\frac{\partial \Phi}{\partial E} \quad (23)$$

Substituting (23) in (21) we get

$$\frac{\partial \Phi}{\partial t} + \frac{\partial \Phi}{\partial E} * \left(\left(\sum_{i=1}^{m}(a_i cos(w_i t) + b_i sin(w_i t))\right)\right) + \frac{1}{4}\left[\frac{\partial \Phi}{\partial E}\right]^2 = 0 \quad (24)$$

According to the (20) let us find the solution of (24) in the form:

$\Phi(t,E) = K(t)E$ (25)

So $\frac{\partial \Phi}{\partial t} = \frac{dK(t)}{dt}E$ (26)

as

$\frac{\partial \Phi}{\partial t} = \frac{dK(t)}{dt}E$ (27)

Substituting (26) and (27) in (24) we can get general form of the equation. If we solve it we get residual error function. It gives us optimal feedback control function:

$$\frac{E(t_1)}{E(t)}\left[\frac{\partial E}{\partial t}\right] + \frac{E(t_1)}{E(t)}\left(\sum_{i=1}^{m}(a_i cos(w_i t) + b_i sin(w_i t))\right) + \frac{1}{4}\frac{E(t_1)}{E(t)} = 0 \quad (28)$$

If $E(t_0)=E_0=0,8777255$ from maximal error range amplitude and round frequencies are equal to 0,07; 1,05;1,48; 1,7; 2,25; 2,60; 3,25 (calculated by "SCAN"); and theirs' amplitude thicknesses are equal to 0,223607; 0,3; 0,3; 0,4472; 0,4472; 0,547; 0,387 from maximal error range amplitude then Runge-Kutta 4th order method gives us the control strategy graph given at fig.1

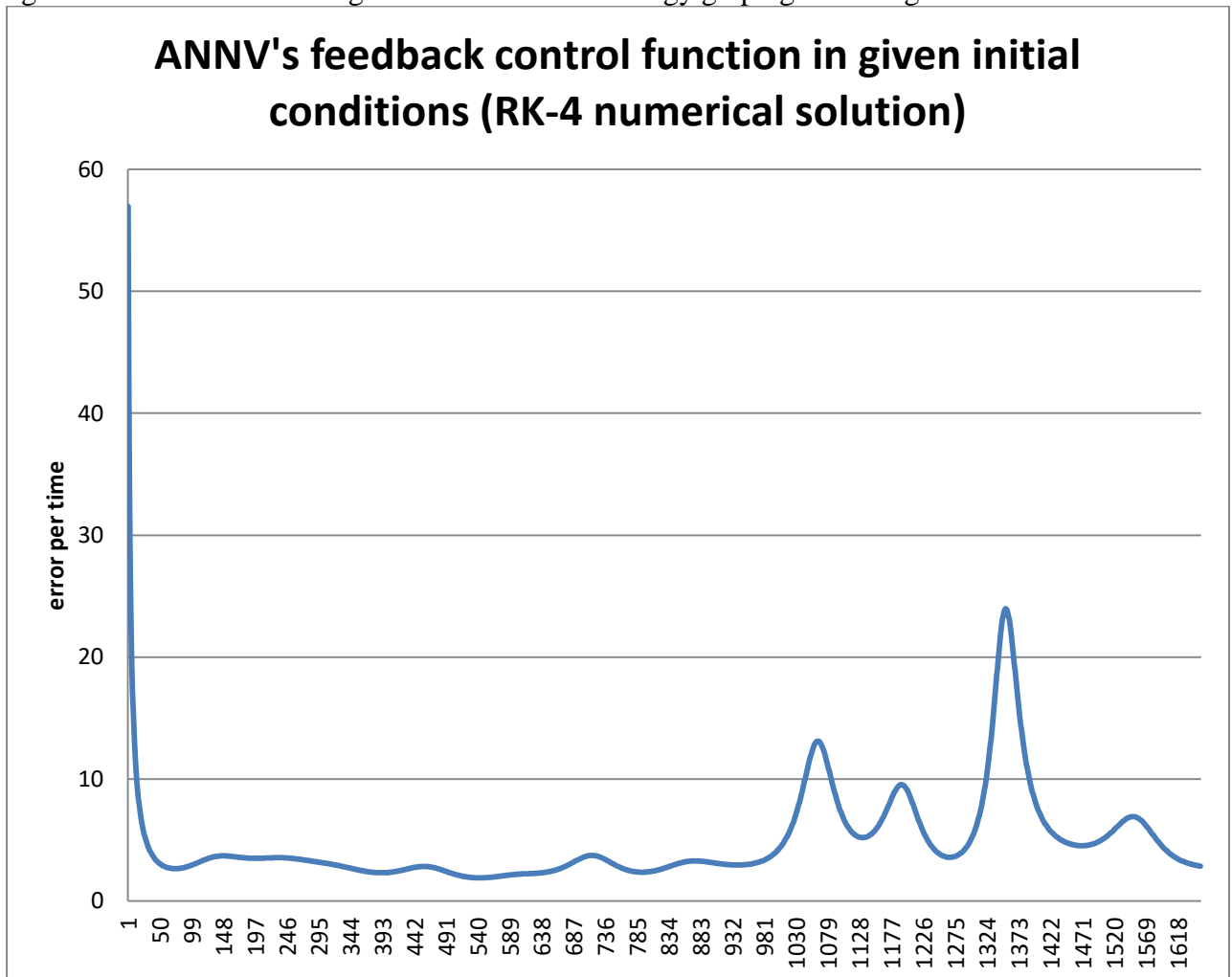

Fig.1 – ANNV control strategy (control function). "X"-axe is time axe.

Conclusions:

1. During the time increasing period the required intensity of control function is to change quasi-periodically (see fig.1)

2. The more harmonics are in the ANNV error evolution equation (4) the more precise the solution is.

3. In general case ANNV is not a fully controllable system, however the solution we have found allows us to consider ANNVs as output-controllable systems.

4. In spite of founded control function it is a very difficult problem to realize the control strategy via algorithms operated with interneuron links changing, weight gradient changing, learning speed by ANNV features and instruments [2]. To realize the control curve it is necessary to conduct experiments to analyze the influence of different weight ensembles and groups on the ANNV error spectrograms. This problem calls for further investigation.

References

1. *Hopfield, J.J.* Neural networks and physical systems with emergent collective computational abilities.//Proc. Natl.Acad.Sci.USA, vol.79 pp.2554-2558, Biophysics, April 1982.
2. *Ignatenkov A.V., Olshansky A.M.* An artificial neural network application to process scheduling in train schedule composing.//Modern IT and IT-education. Vol.2 (11). – 2015. - Moscow State University. – 2015., 614 p., pp.50-55[*Ignatenkov A.V., Olshansky A.M.* Primenenie iskusstvennoi neyronnoi seti dlya postroeniya protsessov na primere grafika dvizheniya poezdov//Sovremennye informazionnie technologii i IT-obrazovanie. T.2(11), 2015.//Moscow, izd-vo VMK MGU, 614 p., p.50-55.]
3. *Fahotimi, O., Dembo, A., Kailath, T.* Neural network weight matrix synthesis using optimal control techniques.//USA, Stanford,1989.//Advances in Neural Information Processing Systems-2 (NIPS-2)//USA, Denver, Colorado.
4. *Becerikli, Y., Konar, A.-F., Samad, T.* Intelligent optimal control with dynamic neural networks.//www.elsevier.com/locate/neunet; Neural Networks 16(2003), pp.251-259.
5. *Tsoy Yu.R.* Neuroevolutionary algorithm and software to image mining. PhD theses. – Tomsk, Tomsky polytechnic university, 2007// - 209 p.[*Tsoy Yu.R.* Neyroevolutzionnyi algoritm I programmnie sredstva dlya obrabotki izobrazhenii: dissertatsia na soiskanie uchenoi stepeni candidata teckchnicheskih nauk. – Tomsk, Tomsky politechnicheskiy universitet, 2007. – 209 p.]
6. *Galimyanov F.A., Gafarov F.M., Khusnutdinov N.R.* A neural network's growth model.//Mathematical modelling. 2011 - #3. – vol.23, p.101-108. [*Galimyanov F.A., Gafarov F.M., Khusnutdinov N.R.* Model rosta neyronnoy seti.//Matematicheskoe modelirovanie. 2011 - #3. – vol.23, p.101-108.]
7. *Andreeva E.A., Pustarnakova Yu. A.* A mathematical model of an artificial neural network with delays.//Program products and systems. – 2001. - #3., pp.6-9[*Andreeva E.A., Pustarnakova Yu. A.* Matematicheskaya model neyronnoy seti s zapazdivaniem//Programmnye produkty i sistemy. – 2001. - #3., pp.6-9]
8. *Haykin S.*, Neural networks//Moscow, Williams publishing, 2nd edition. 2006., 1104 p. – ISBN 5-8459-0890-6 (rus).
9. *Olshansky A.M.* Some features of artificial network control problems foundation.//Modern principles, methods and automation of transport control systems. Proceedings of the International scientific and practical conference - Moscow, Moscow state railway university, 2016. – 192 p., pp.111-114. [*Olshanksy A.M.* Nekotorie osobennosti postanovki zadach upravlenija neyronnimi setyami//Sovremennye metody, principy i sistemy avtomatizatsii upravlenija na transporte: sbornik trudov Mezhdunarodnoi nauchno-prakticheskoi konferencii (19-20 april 2016, Nizhniy Novgorod, - Moscow, Moscow state rail university, 2016 - 192 p., pp.111-114]
10. *Pantellev A.V., Bortakovsky A.S.* Control theory in examples and problems. Moscow, Higher Education, 2003. - 583 p., p.487. [*Pantellev A.V., Bortakovsky A.S.* Teoriya upravlenija v primerakh i zadachakh. Moscow, Viscshaya shkola, 2003. – 583 p., p.487]
11. *Tanaev V.S., Shkurba V.V.* Introduction to schedule theory. - Moscow, The Science, 1975. [*Tanaev V.S., Shkurba V.V.* Vvedenie v teoriyu raspisanij. – Moscow, Nauka, 1975.]